\documentstyle[12pt,leqno]{article}
\font\erlss=eufm9
\newcommand{\Lss}{\mbox{\erlss L}}
\font\erl=eufm10 at 12pt
\newcommand{\LL}{\mbox{\erl L}}
\font\ers=eusm10 at 12pt
\newcommand{\HH}{\mbox{\ers H}}
\font\erss=eusm10
\newcommand{\HHs}{\mbox{\erss H}}
\font\ms=msam8 at 10pt
\newcommand{\eql}{\mbox{\thinspace\thinspace
{\ms5}\thinspace\thinspace}}
\newcommand{\eqg}{\mbox{\thinspace\thinspace
{\ms=}\thinspace\thinspace}}
\font\mss=msam8
\newcommand{\eqls}{\mbox{{\mss 5}}}

\font\msb=msbm8 at 12pt
\newcommand{\R}{\mbox{\msb R}}
\newcommand{\C}{\mbox{\msb C}}
\newcommand{\N}{\mbox{\msb N}}
\begin{document}
\vspace*{15mm}
\begin{center}
{\large \bf ON THE {\boldmath$H$}-FUNCTION}\\[5mm]
{\large Anatoly A. Kilbas}\\
{\it Department of Mathematics and Mechanics$,$ Belarusian State
University\\
Minsk 220050$,$ Belarus}\\[5mm]
{\large Megumi Saigo}\\
{\it Department of Applied Mathematics$,$ Fukuoka University\\
Fukuoka 814-0180$,$ Japan}
\end{center}
\vspace*{15mm}
\begin{abstract}{The paper is devoted to study the $H$-function
defined by the 
Mellin-Barnes integral 
$$H^{m,n}_{\thinspace p,q}(z)={\frac1{2\pi i}}\int_{\Lss}
\HHs^{m,n}_{\thinspace p,q}(s)z^{-s}ds,$$
where the function $\HH^{m,n}_{\thinspace p,q}(s)$ is a certain
ratio of 
products of Gamma functions with the argument $s$ and the contour
$\LL$ is 
specially chosen. The conditions for the existence of 
$H^{m,n}_{\thinspace p,q}(z)$ are discussed and explicit power
and 
power-logarithmic series expansions of $H^{m,n}_{p,q}(z)$ near
zero and 
infinity are given. The obtained results define more precisely
the known 
results.}
\\\par
{\bf Key words}: $H$-function\par
{\bf AMS(MOS) Subject Classification}: 33C40
\end{abstract}
\vspace{7mm}
\begin{flushleft}
{\bf 1. \ Introduction}
\end{flushleft}
\setcounter{section}{1}
\setcounter{equation}{0}
\par
This paper deals with the $H$-function $H^{m,n}_{p,q}(z)$
introduced by 
Pincherle in 1888 (see [3, Section\thinspace1.19]). Interest in 
this function appeared in 1961, when Fox [4] investigated such a
function as 
symmetrical Fourier kernel. Therefore, the $H$-function is 
often called Fox's $H$-function. For integers $m,n,p,q$ such that

$0\eql m\eql q, 0\eql n\eql p$ and $a_i, b_j \in{\C}$ and 
$\alpha_i, \beta_j \in{\R}_+=(0,\infty)$ $(1\eql i\eql p,$ $1\eql
j\eql q)$ 
the function is defined by the Mellin-Barnes integral
\begin{eqnarray}
H^{m,n}_{p,q}(z)={\frac1{2\pi
i}}\int_{\Lss}\HH^{m,n}_{p,q}(s)z^{-s}ds,
\end{eqnarray}
where
\begin{eqnarray}
&&\HH^{m,n}_{p,q}(s)={\frac{\displaystyle{\prod^m_{j=1}\Gamma(b_j
+\beta_js)
\prod^n_{i=1}\Gamma(1-a_i-\alpha_is)}}
{\displaystyle{\prod^p_{i=n+1}\Gamma(a_i+\alpha_is)\prod^q_
{j=m+1}
\Gamma(1-b_j-\beta_js)}}},
\end{eqnarray}
the contour $\LL$ is specially chosen and an empty product, if it
occurs, is 
taken to be one. The theory of this function may be found in
[10], [2], [1], 
[9, Chapter 2], [12, Chapter 1] and [11, Section\thinspace8.3].
We only 
indicate that most of the elementary and special functions are
particular 
cases of the $H$-function $H^{m,n}_{p,q}(z)$. In particular, if
$\alpha$'s 
and $\beta$'s are equal to 1, the $H$-function (1.1) reduces to
Meijer's 
$G$-function $G^{m,n}_{p,q}(z)$.   
\par
The conditions of the existence of the $H$-function can be made
by inspecting 
the convergence of the integral (1.1), which depend on 
the selection of the contour $\LL$ and the relations between
parameters 
$a_i,\alpha_{i}\ (i=1,\cdots,p)$ and 
$b_j,\beta_j\ (j=1,\cdots,q)$. Especially, the relations may
depend on the 
numbers $\Delta,\delta$ and $\mu$ defined by  
\begin{eqnarray}
\Delta&\hspace{-2.5mm}=&\hspace{-2.5mm}\sum^{q}_{j=1}\beta_{j}-
\sum^{p}_{i=1}
\alpha_{i},\\[2mm]
\delta&\hspace{-2.5mm}=&\hspace{-2.5mm}\prod^{p}_{i=1}\alpha_{i}^
{-\alpha_{i}}
\prod^{q}_{j=1}\beta_{j}^{\beta_{j}}.\\[2mm]
\mu&\hspace{-2.5mm}=&\hspace{-2.5mm}\sum^{q}_{j=1}b_{j}-\sum^{p}_
{i=1}a_{i}
+{\frac{p-q}{2}}.
\end{eqnarray}
Such a selection of $\LL$ and the relations on parameters are
indicated in 
the handbook [11, Section\thinspace8.3.1], but some of the
results there 
needs correction. In this paper we would like to give such a
correction in 
the following cases: 
\par
{\bf (a)} \ $\Delta \eqg 0$ and the contour $\LL=\LL_{-\infty}$
in (1.1) runs 
from $-\infty +i\varphi_{1}$ to $-\infty +i\varphi_{2},\
\varphi_{1}
<\varphi_{2},$ such that the poles of $\Gamma(b_j+\beta_js)\
(j=1,\cdots,m)$ 
lie on the left of $\LL_{-\infty}$ and those of
$\Gamma(1-a_i-\alpha_is)\ 
(i=1,\cdots,n)$ on the right of $\LL_{-\infty}$. 
\par
{\bf (b)} \ $\Delta \eql 0$ and the contour $\LL=\LL_{+\infty}$
in (1.1) runs 
from $+\infty +i\varphi_{1}$ to $+\infty +i\varphi_{2},\
\varphi_{1}
<\varphi_{2},$ such that the poles of $\Gamma(b_j+\beta_js)\
(j=1,\cdots,m)$ 
lie on the left of $\LL_{+\infty}$ and thoes of
$\Gamma(1-a_i-\alpha_is)\ 
(i=1,\cdots,n)$ on the right of $\LL_{+\infty}$.
\par
Our results are based on the asymptotic behavior of the function 
$\HH^{m,n}_{p,q}(s)$ given in (1.2) at infinity. Using the
behavior 
and following [1], we give the series representation of
$H^{m,n}_{p,q}(z)$ 
via residues of the integrand $\HH^{m,n}_{p,q}(s)z^{-s}$. In this
way we 
simplify the proof of Theorem 1 in [1] by applying the former
results to 
find the explicit series expansions of $H^{m,n}_{p,q}(z)$. Such
power 
expansions, as corollaries of the results from [1], were
indicated in 
[12, Chapter 2.2] (see also [11, Section\thinspace8.3.1]),
provided that 
the poles of Gamma-functions $\Gamma(b_j+\beta_js)\
(j=1,\cdots,m)$ and 
$\Gamma(1-a_i-\alpha_is)\ (i=1,\cdots,n)$ do not coincide
\begin{eqnarray}
&&\beta_{j}(a_{i}-1-k)\neq
\alpha_{i}(b_{j}+l)\quad(i=1,\cdots,n;\ 
j=1,\cdots ,m;\ k,l\in\N_0=\{0,1,2,\cdots\})
\end{eqnarray}
in the cases:
\par
{\bf (c)} \ $\Delta >0$ with $z\neq 0$ or $\Delta=0$ with
$0<|z|<\delta$, and 
the poles of Gamma-functions 
$\Gamma(b_j+\beta_js)\ (j=1,\cdots,m)$ are simple: 
\begin{eqnarray}
&&\beta_{j}(b_{i}+k)\neq \beta_{i}(b_{j}+l),\quad(i\neq j,\
i,j=1,\cdots,m;\ 
k,l\in\N_0);
\end{eqnarray}
\par
{\bf (d)} \ $\Delta <0$ with $z\neq 0$ or $\Delta=0$ with
$|z|>\delta$, and 
the poles of Gamma-functions $\Gamma(1-a_i-\alpha_is)\
(i=1,\cdots,n)$ are 
simple:
\begin{eqnarray}
&&\alpha_{j}(1-a_{i}+k)\neq \alpha_{i}(1-a_{j}+l),\quad(i\neq j,\
i,j=1,\cdots,n;\ k,l\in\N_0).
\end{eqnarray}
\par
When the poles of Gamma-functions in (c) and (d) coincide,
explicit series 
expansions of $H^{m,n}_{p,q}(z)$ should be more complicated
power-logarithmic 
expansions. Such expansions in particular cases of the Meijer's
$G$-functions 
$G^{p,0}_{0,p}$ and $G^{p,0}_{p,p}$ and of the $H$-functions
$H^{p,0}_{0,p}$ 
and $H^{p,0}_{p,p}$ were given in [7] and [8], respectively. 
\par
We obtain the explicit expansions of the $H$-function of general
form 
$H^{m,n}_{p,q}(z)$ under the conditions in (1.6). We show that, 
if the poles of the Gamma-functions $\Gamma(b_j+\beta_js)\
(j=1,\cdots,m)$ 
and $\Gamma(1-a_i-\alpha_is)\ (i=1,\cdots,n)$ coincide in the
cases (c) and 
(d), respectively, then the $H$-function (1.1) has
power-logarithmic series 
expansions. In particular, we give the asymptotic expansions of 
$H^{m,n}_{p,q}(z)$ near zero. We note that the obtained results
will be 
different in the cases when either $\Delta \eqg 0$ or $\Delta
\eql 0$.
\par
The paper is organized as follows. Section 2 is devoted to the
conditions of 
the existence of the $H$-function (1.1) which are based on the
asymptotic 
behavior of $\HH^{m,n}_{p,q}(s)$ at infinity. Here we also give
the 
representations of (1.1) via the residues of the integrand. The
latter result 
is applied in Sections 3 and 4 to obtain the explicit power and 
power-logarithmic series expansions of $H^{m,n}_{p,q}(z)$ and, in
particular, 
its asymptotic estimates near zero.
\\\\
%
%
\begin{flushleft}
{\bf 2. Existence and Representations of $H^{m,n}_{p,q}(z)$}
\end{flushleft}
\setcounter{section}{2}
\setcounter{equation}{0}
\par
First we give the asymptotic estimate of Gamma function $\Gamma
(z),\ z=x+iy,$ 
[3, Chapter 1] at infinity on lines parallel to the coordinate
axes.
\\\par
{\bf Lemma 1.} \ {\sl Let $z=x+iy$ with $x,y\in{\R}=(-\infty,
\infty)$. 
Then the following asymptotic estimates at infinity are valid$:$
\begin{eqnarray}
&&|\Gamma (x+iy)|\sim\sqrt{2\pi}|x|^{x-1/2}e^{-x-\pi(1-{\rm
\scriptstyle sign}
\thinspace x)y/2}\quad(|x|\to\infty;\ y\ne0\mbox{ if }x<0)
\end{eqnarray}
and}
\begin{eqnarray}
&&|\Gamma
(x+iy)|\sim\sqrt{2\pi}|y|^{x-1/2}e^{-x-\pi|y|/2}\quad(|y|\to
\infty).
\end{eqnarray}
\par
{\bf Proof.} Applying the Stirling formula [3, 1.18(2)]
\begin{eqnarray}
&&\Gamma (z)\sim\sqrt{2\pi}e^{(z-1/2)\log
z}e^{-z}\quad(|z|\to\infty;\ 
|\arg(z)|<\pi),
\end{eqnarray}
we have
\begin{eqnarray}
|\Gamma (x+iy)|&\hspace{-2.5mm}\sim&\hspace{-2.5mm}\sqrt{2\pi}
\left|e^{(x+iy-1/2)[\log |x+iy|+i\arg
(x+iy)]}e^{-(x+iy)}\right|\\[2mm]
&\hspace{-2.5mm}\sim&\hspace{-2.5mm}\sqrt{2\pi}|x+iy|^{x-1/2}e^{-
x-y
\arg (x+iy)}\quad(|x+iy|\to \infty;\ y\ne0\mbox{ if
}x\eql0).\nonumber
\end{eqnarray}
\par
Let $y\in{\R}$ be fixed and $|x|\to\infty$. Then $|x+iy|\sim|x|$,
and $\arg
(x+iy)\to 0$ as $x\to+\infty$ and $\arg(x+iy)\to\pi$ as 
$x\to -\infty$. Therefore, (2.4) implies
\begin{eqnarray}
&&|\Gamma
(x+iy)|\sim\sqrt{2\pi}|x|^{x-1/2}e^{-x}\quad(x\to+\infty)
\end{eqnarray}
and
\begin{eqnarray}
&&|\Gamma (x+iy)|\sim\sqrt{2\pi}|x|^{x-1/2}e^{-x-\pi
y}\quad(x\to-\infty;\ y
\ne0),
\end{eqnarray}
which yield (2.1).
\par
Turning to the case $x\in{\R}$ being fixed and $|y|\to\infty$, we
find 
$|x+iy|\sim |y|$ and $\arg (x+iy)\to \pi/2$ as $y\to \infty$ and 
$\arg (x+iy)\to -\pi/2$ as $y\to -\infty$. Thus (2.4) implies
(2.2).
\\\par
{\bf Remark 1.} \ The relation [3, (1.18.6)] needs correction
with addition 
of the multiplier $e^{x}$ in the left hand side and it must be
replaced by
\begin{eqnarray}
&&\lim_{|y|\to \infty}|\Gamma (x+iy)|e^{x+\pi
|y|/2}|y|^{1/2-x}=\sqrt{2\pi}.
\end{eqnarray}
\\\par
Next assertion gives the asymptotic behavior of
$\HH^{m,n}_{p,q}(s)$ defined 
in (1.2) at infinity on lines parallel to the real axis.
\\\par
{\bf Lemma 2.} \ {\sl Let $\Delta,\delta$ and $\mu $ be given by
$(1.3)$ to 
$(1.5)$ and let $t,\sigma \in {\R}$. Then there hold the
estimates
\begin{eqnarray}
&&|\HH^{m,n}_{p,q}(t+i\sigma )|\sim
A\left({\frac{e}{t}}\right)^{-\Delta t}
\delta^{t}t^{{\rm Re}(\mu)}\quad(t\to+\infty)
\end{eqnarray}
with
\begin{eqnarray}
&&A=(2\pi)^{m+n-(p+q)/2}e^{q-m-n}\ 
{\frac{\displaystyle{\prod^{q}_{j=1}\left[(\beta_{j})^{{\rm
Re}(b_{j})-1/2}
e^{-{\rm Re}(b_{j})}
\right]}}{\displaystyle{\prod^{p}_{i=1}\left[(\alpha_{i})^{{\rm
Re}(a_{i})-1/2}
e^{-{\rm Re}(a_{i})}\right]}}}
{\frac{\displaystyle{\prod^{n}_{i=1}e^{\pi [\sigma \alpha_{i}
+{\rm Im}(a_{i})]}}}{\displaystyle{\prod^{q}_{j=m+1}e^{\pi
[\sigma \beta_{j}+{\rm Im}(b_{j})]}}}},
\end{eqnarray}
and 
\begin{eqnarray}
&&|\HH^{m,n}_{p,q}(t+i\sigma )|\sim
B\left({\frac{e}{|t|}}\right)^{\Delta |t|}
\delta^{-|t|}|t|^{{\rm Re}(\mu)}
\quad(t\to-\infty)
\end{eqnarray}
with}
\begin{eqnarray}
&&B=(2\pi)^{m+n-(p+q)/2}e^{q-m-n}\ {\frac{\displaystyle
{\prod^{q}_{j=1}\left[(\beta_{j})^{{\rm Re}(b_{j})-1/2}e^{-{\rm
Re}(b_{j})}
\right]}}{\displaystyle
{\prod^{p}_{i=1}\left[(\alpha_{i})^{{\rm Re}(a_{i})-1/2}e^{-{\rm
Re}(a_{i})}
\right]}}}
{\frac{\displaystyle{\prod^{p}_{i=n+1}e^{\pi [\sigma \alpha_{i}
+{\rm Im}(a_{i})]}}}{\displaystyle{\prod^{m}_{j=1}
e^{\pi [\sigma \beta_{j}+{\rm Im}(b_{j})]}}}}.
\end{eqnarray}
\par
{\bf Proof.} \ By virtue of (2.1), we have, for a variable
$s=t+i\sigma$ and 
a complex constant $k=c+id$,
\begin{eqnarray}
&&|\Gamma
(s+k)|\sim\sqrt{2\pi}t^{t+c-1/2}e^{-(t+c)}\quad(t\to+\infty)
\end{eqnarray}
and
\begin{eqnarray}
&&|\Gamma
(s+k)|\sim\sqrt{2\pi}|t|^{t+c-1/2}e^{-(t+c)}e^{-\pi(\sigma+d)}
\quad(t\to-\infty).
\end{eqnarray}
Substituting these estimates into (1.2) and using (1.3) to (1.5),
we obtain 
(2.8) and (2.10).
\\\par
{\bf Remark 2.} \ The asymptotic estimate of the function
$\HH^{m,n}_{p,q}(s)$ 
at infinity on lines parallel to the imaginary axis 
$\HH^{m,n}_{p,q}(\sigma +it)$ as $|t|\to \infty$ was given in our
paper with 
Shlapakov [5].
\\\par
By appealing to Lemma 2, we give conditions of the existence of
the 
$H$-function (1.1) with the contour $\LL$ being chosen as
indicated in (a) 
and (b) in Section 1.  
\\\par
{\bf Theorem 1.} \ {\sl Let $\Delta,\delta$ and $\mu $ be given
by $(1.3)$ 
to $(1.5)$. Then the function $H^{m,n}_{p,q}(z)$ defined by
$(1.1)$ and 
$(1.2)$ exists in the following cases}:
\begin{eqnarray}
&&\LL=\LL_{-\infty},\quad\Delta>0,\quad z\neq
0;\label{1.2.14}\\[2mm]
&&\LL=\LL_{-\infty},\quad\Delta=0,\quad0<\left|z\right|<\delta;\l
abel{1.2.15}
\\[2mm]
&&\LL=\LL_{-\infty},\quad\Delta=0,\quad\left|z\right|=\delta,
\quad
{\rm Re}(\mu)<-1;\label{1.2.16}\\[2mm]
&&\LL=\LL_{+\infty},\quad\Delta<0,\quad z\neq
0;\label{1.2.17}\\[2mm]
&&\LL=\LL_{+\infty},\quad\Delta=0,\quad\left|z\right|>\delta;
\label{1.2.18}\\[2mm]
&&\LL=\LL_{+\infty},\quad\Delta=0,\quad\left|z\right|=\delta,
\quad{\rm Re}(\mu)<-1.
\end{eqnarray}
\par
{\bf Proof.} Let us first consider the case (a) for which $\Delta
\eqg 0$ and 
$\LL=\LL_{-\infty}$. We have to investigate the convergence of
the integral 
(1.1) on the lines 
\begin{eqnarray}
&&l_{1}=\{t\in{\R}:t+i\varphi_1\}\quad\mbox{and}\quad
l_{2}=\{t\in{\R}:
t+i\varphi_2\}\quad\mbox{for}\quad\varphi_1
<\varphi_2
\end{eqnarray}
as $t\to-\infty$. According to (2.10), we have the following
asymptotic 
estimate for the integrand of (1.1):
\begin{eqnarray}
&&|{\HH}^{m,n}_{p,q}(s)z^{-s}|\sim B_{i}e^{\varphi_i\arg
z}\left({\frac{e}
{|t|}}\right)^{\Delta |t|}\left({\frac{|z|}{\delta}}\right)^{|t|}
|t|^{{\rm Re}(\mu)}\quad(t\to-\infty;\ t\in l_{i}\ (i=1,2)),
\end{eqnarray}
where $B_{1}$ and $B_{2}$ are given by (2.11) with $\sigma$ being
replaced by 
$\varphi_{1}$ and $\varphi_{2}$, respectively. It follows from
(2.21) that
the integral (1.1) is convergent if and only if one of the
conditions in
(2.14) to (2.16) is satisfied.
\par
In the case (b), $\Delta \eql 0$ and the contour $\LL$ is taken
to be 
$\LL_{+\infty}$. Then we have to investigate the convergence 
of the integral (1.1) on the lines $l_1$ and $l_2$ in (2.20), as 
$t\to +\infty$. By virtue of (2.8) and (2.9), we have the
asymptotic 
estimate:
\begin{eqnarray}
&&|{\HH}^{m,n}_{p,q}(s)z^{-s}|\sim A_{i}e^{\varphi_i\arg
z}\left({\frac{e}{t}}
\right)^{-\Delta t}\left({\frac{\delta}{|z|}}\right)^{t}t^{{\rm
Re}(\mu)}
\quad(t\to +\infty;\ t\in l_{i}\ (i=1,2)),
\end{eqnarray}
where $A_{1}$ and $A_{2}$ are given by (2.9) with $\sigma$ being
replaced by 
$\varphi_{1}$ and $\varphi_{2}$, respectively. Thus (2.22)
implies that the 
integral (1.1) converges if and only if one of the conditions in
(2.17) to 
(2.19) holds.
\\\par
{\bf Corollary 1.} \ {\sl The estimate $(2.21)$ holds for $t\to
-\infty$ 
uniformly on sets which have a positive distance to the points   
\begin{eqnarray}
&&b_{jl}=-{\frac{b_{j}+l}{\beta_{j}}}\quad(j=1,\cdots,m;\
l\in\N_0).
\end{eqnarray}
and do not contain points to the right of $\LL_{-\infty}$.
\par
The estimate $(2.22)$ holds for $t\to +\infty$ uniformly on sets
which have 
a positive distance to the points   
\begin{eqnarray}
&&a_{ik}={\frac{1-a_{i}+k}{\alpha_{i}}}\quad(i=1,\cdots,n;\
k\in\N_0).
\end{eqnarray}
and do not contain points to the left of} $\LL_{+\infty}$.
\\\par
{\bf Remark 3.} \ The conditions for the existence of the
$H$-function (1.1)  
\begin{eqnarray}
&&\sum^n_{i=1}\alpha_i-\sum^p_{i=n+1}\alpha_i+\sum^m_{j=1}\beta_j
-\sum^q_{j=m+1}\beta_j\eqg0,\quad{\rm Re}(\mu)<0
\end{eqnarray}
given in [11, Section\thinspace8.3.1] in the cases when
$\LL=\LL_{-\infty},
\Delta =0,|z|=\delta$ ([8, Section\thinspace8.3.1.3)]) and
$\LL=\LL_{+\infty},
\Delta =0,|z|=\delta$ ([8, Section\thinspace8.3.1.4)]) can be
replaced by the 
condition
\begin{eqnarray}
&&{\rm Re}(\mu)<-1.
\end{eqnarray}
\\\par
The following statement follows from Theorem 1, Corollary 1 and
the theory of 
residues.
\\\par
{\bf Theorem 2.} \ {\bf (A)} \ {\sl If the conditions in $(1.6),$
and $(2.14)$ 
or $(2.15)$ are satisfied$,$ then the $H$-function $(1.1)$ is an
analytic 
function of $z$ in the corresponding domain indicated in $(2.14)$
or $(2.15),$ 
and
\begin{eqnarray}
&&H^{m,n}_{p,q}(z)=\sum^{m}_{j=1}\sum^{\infty}_{l=0}\mathop{\rm
Res}_{s=b_{jl}}
[{\HH}^{m,n}_{p,q}(s)z^{-s}],
\end{eqnarray}
where $b_{jl}$ are given in $(2.23).$
\par
{\bf (B)} \ If the conditions in $(1.6),$ and $(2.17)$ or
$(2.18)$ are 
satisfied$,$ then the $H$-function $(1.1)$ is an analytic
function of $z$ in 
the corresponding domain indicated in $(2.17)$ or $(2.18),$ and 
\begin{eqnarray}
&&H^{m,n}_{p,q}(z)=-\sum^{n}_{i=1}\sum^{\infty}_{k=0}
\mathop{\rm Res}_{s=a_{ik}}[{\HH}^{m,n}_{p,q}(s)z^{-s}],
\end{eqnarray}
where $a_{ik}$ are given in} (2.24).
\\\par
{\bf Remark 4.} \ The first assertion of Theorem 2 was proved in
[1, p.278, 
Theorem 1] for the $H$-function represented by the integral
obtained from 
(1.1) and (1.2) after replacing $s$ by $-s$. The proof of Theorem
1 in [1] is 
complicated and based on Lemma 2 there in which the asymptotic
estimate at 
infinity of the functions $h_{0}(s)$ defined by
\begin{eqnarray}
&&h_{0}(s)={\frac{\displaystyle{\prod^{p}_{i=1}\Gamma
(1-a_{i}+\alpha_{i}s)}}
{\displaystyle{\prod^{q}_{j=1}\Gamma (1-b_{j}+\beta_{j}s)}}}
\end{eqnarray}
is given. But our proof of Theorem 2 along the ideas of [1] is
more simple 
and is based on the asymptotic estimate of ${\HH}^{m,n}_{p,q}(s)$
at infinity 
given in Lemma 2.
\\\\
\begin{flushleft}
{\bf 3. Explicit Power Series Expansions}
\end{flushleft}
\setcounter{section}{3}
\setcounter{equation}{0}
\par
In this section we apply Theorem 2 to obtain explicit power
series expansions 
of the $H$-function (1.1) under the condition (1.6) in the case
of (1.7) or 
(1.8).
\par
First we consider the former case. By Theorem 2(A), we have to
evaluate the 
residues of ${\HH}(s)z^{-s}$ at the points $s=b_{jl}$ given in
(2.23), where 
and in what follows we simplify ${\HH}^{m,n}_{p,q}(s)$ by
${\HH}(s)$. To 
evaluate these residues we use the property of the Gamma-function
[6, (3.30)], 
that is, in a neighbourhood of the poles $z=-k\ (k\in\N_0)$ the
Gamma-function 
$\Gamma (z)$ can be expanded in powers of $z+k=\epsilon$
\begin{eqnarray}
&&\Gamma (z)={\frac{(-1)^{k}}{k!\epsilon}}[1+\epsilon\psi(1+k)
+O(\epsilon^{2})],\quad\mbox{where}\quad\psi
(z)={\frac{\Gamma'(z)}
{\Gamma (z)}}.
\end{eqnarray}
\par
Since the poles $b_{jl}$ are simple, i.e., the conditions in
(1.7) hold,
\begin{eqnarray}
&&\mathop{\rm
Res}_{s=b_{jl}}[{\HH}(s)z^{-s}]=h^{\ast}_{jl}z^{-b_{jl}}
\ (j=1,\cdots,m;\ l\in\N_0),
\end{eqnarray}
where
\begin{eqnarray}
h^{\ast}_{jl}&\hspace{-2.5mm}=&\hspace{-2.5mm}\lim_{s\to b_{jl}}
\left[(s-b_{jl}){\HH}(s)\right]\\[2mm]
&\hspace{-2.5mm}=&\hspace{-2.5mm}{\frac{(-1)^{l}}{l!\beta_{j}}}
{\frac{\displaystyle{\prod^m_{i=1,i\neq j}\Gamma\left(b_i-[b_j+l]
\frac{\beta_i}{\beta_j}\right)\prod^n_{i=1}\Gamma\left(1-a_i+[b_j
+l]
\frac{\alpha_i}{\beta_j}\right)}}{\displaystyle{\prod^p_{i=n+1}
\Gamma\left(a_i-[b_j+l]\frac{\alpha_i}{\beta_j}\right)\prod^q_{i=
m+1}
\Gamma\left(1-b_i+[b_j+l]\frac{\beta_{i}}{\beta_j}\right)}}}.
\nonumber
\end{eqnarray}
\par
Thus we obtain
\\\par
{\bf Theorem 3.} \ {\sl Let the conditions in $(1.6)$ and $(1.7)$
be satisfied 
and let either $\Delta>0,z\neq0$ or $\Delta=0,0<|z|<\delta$. Then
the 
$H$-function $(1.1)$ has the power series expansion   
\begin{eqnarray}
&&H^{m,n}_{p,q}(z)=\sum^{m}_{j=1}\sum^{\infty}_{l=0}h^{\ast}_{jl}
z^{(b_{j}+l)/
\beta_{j}},
\end{eqnarray}
where the constants $h^{\ast}_{jl}$ are given by} (3.3).
\\\par
{\bf Corollary 2.} \ {\sl If the conditions in $(1.6)$ and
$(1.7)$ are 
satisfied and $\Delta\eqg0,$ then $(3.4)$ gives the asymptotic
expansion of 
$H^{m,n}_{p,q}(z)$ near zero and the main terms of this
asymptotic formula 
have the form$:$
\begin{eqnarray}
&&H^{m,n}_{p,q}(z)=\sum^m_{j=1}\left[h^*_jz^{b_j/\beta_j}+O\left(
z^{(b_j+1)/
\beta_j}\right)\right]\quad(z\to0),
\end{eqnarray}
where}
\begin{eqnarray}
&&h^*_j\equiv
h^*_{j0}={\frac1{\beta_j}}{\frac{\displaystyle{\prod^m_{i=1,i
\neq j}\Gamma\left(b_i-\frac{b_j\beta_i}{\beta_j}\right)
\prod^n_{i=1}\Gamma\left(1-a_i+\frac{b_j\alpha_i}{\beta_j}\right)
}}
{\displaystyle{\prod^p_{i=n+1}\Gamma\left(a_i-\frac{b_j\alpha_i}
{\beta_j}\right)\prod^q_{i=m+1}\Gamma\left(1-b_i+\frac{b_j\beta_i
}{\beta_j}
\right)}}}.
\end{eqnarray}
\\\par
{\bf Corollary 3.} \ {\sl Let the conditions in $(1.6)$ and
$(1.7)$ be 
satisfied$,$ and let $\Delta\eqg0$ and $j_0\ (1\eql j_0\eql m)$
be an integer 
such that
\begin{eqnarray}
&&{\frac{{\rm Re}(b_{j_0})}{\beta_{j_0}}}=\min_{1\eqls j\eqls m}
\left[{\frac{{\rm Re}(b_j)}{\beta_j}}\right].
\end{eqnarray}
Then there holds the asymptotic estimate$:$
\begin{eqnarray}
&&H^{m,n}_{p,q}(z)=h^*_{j_0}z^{b_{j_0}/\beta_{j_0}}+o\left(z^{b_{
j_0}/
\beta_{j_0}}\right)\quad(z\to0),
\end{eqnarray}
where $h^*_{j_0}$ is given by $(3.6)$ with $j=j_0$. In
particular}, 
\begin{eqnarray}
&&H^{m,n}_{p,q}(z)=O(z^\rho)\quad(z\to0)\quad\mbox{with}\quad\rho
=\min_{1\eqls j\eqls m}
\left[{\frac{{\rm Re}(b_j)}{\beta_j}}\right].
\end{eqnarray}
\\\par
Now we consider the case (1.8) when the poles of the
Gamma-functions 
$\Gamma(1-a_i-\alpha_is)\ (i=1,\cdots,n)$ are simple. By (3.1),
evaluating 
the residues of ${\HH}(s)z^{-s}$ at the points $a_{ik}$ given in
(2.24) we 
have similarly to the previous argument that
\begin{eqnarray}
&&\mathop{\rm Res}_{s=a_{ik}}[{\HH}(s)z^{-s}]=-h_{ik}z^{-a_{ik}}
\quad(i=1,\cdots,n;\ k\in\N_0),
\end{eqnarray}
where $a_{ik}$ are given by (2.24) and 
\setcounter{equation}{10}
\begin{eqnarray}
h_{ik}&\hspace{-2.5mm}=&\hspace{-2.5mm}\lim_{s\to
a_{ik}}\left[-(s-a_{ik})
{\HH}(s)\right]\\[2mm]
&\hspace{-2.5mm}=&\hspace{-2.5mm}{\frac{(-1)^k}{k!\alpha_i}}
{\frac{\displaystyle{\prod^m_{j=1}\Gamma\left(b_j+[1-a_i+k]\frac{
\beta_j}
{\alpha_i}\right)
\prod^n_{j=1,j\neq
i}\Gamma\left(1-a_j-[1-a_{i}+k]\frac{\alpha_j}{\alpha_i}
\right)}}
{\displaystyle{\prod^p_{j=n+1}\Gamma\left(a_j+[1-a_{i}+k]\frac{\a
lpha_j}
{\alpha_i}\right)
\prod^q_{j=m+1}\Gamma\left(1-b_j-[1-a_{i}+k]\frac{\beta_j}{\alpha
_i}\right)}}}.
\nonumber
\end{eqnarray}
\par
Thus from Theorem 2(B) we have
\\\par
{\bf Theorem 4.} \ {\sl Let the conditions in $(1.6)$ and $(1.8)$
be satisfied 
and let either $\Delta<0,z\neq0$ or $\Delta=0,|z|>\delta$. Then
the 
$H$-function $(1.1)$ has the power series expansion   
\begin{eqnarray}
&&H^{m,n}_{p,q}(z)=\sum^{n}_{i=1}\sum^{\infty}_{k=0}h_{ik}z^{(a_{
i}-k-1)/
\alpha_{i}},
\end{eqnarray}
where the constants $h_{ik}$ are given by} (3.11).
\\\par
{\bf Corollary 4.} \ {\sl If the conditions in $(1.6)$ and
$(1.8)$ are 
satisfied and $\Delta<0,$ then $(3.12)$ gives the asymptotic
expansion of 
$H^{m,n}_{p,q}(z)$ near infnity and the main terms of this
asymptotic formula 
have the form$:$
\begin{eqnarray}
&&H^{m,n}_{p,q}(z)=\sum^n_{i=1}\left[h_iz^{(a_i-1)/\alpha_i}+O
\left(z^{(a_i-2)/
\alpha_i}\right)\right]\quad(|z|\to\infty),
\end{eqnarray}
where
\begin{eqnarray}
&&h_i\equiv h_{i0}=\frac{1}{\alpha_i}
{\frac{\displaystyle{\prod^m_{j=1}\Gamma\left(b_j-[a_i-1]\frac
{\beta_j}
{\alpha_i}\right)
\prod^n_{j=1,j
\neq i}\Gamma\left(1-a_j+[a_i-1]\frac{\alpha_j}{\alpha_i}
\right)}}
{\displaystyle{\prod^p_{j=n+1}\Gamma\left(a_j-[a_i-1]\frac{\alpha
_j}{\alpha_i}
\right)
\prod^q_{j=m+1}\Gamma\left(1-b_j+[a_i-1]\frac{\beta_j}{\alpha_i}
\right)}}}.
\end{eqnarray}
\\\par
{\bf Corollary 5.} \ {\sl Let the conditions in $(1.6)$ and
$(1.8)$ be 
satisfied$,$ and let $\Delta<0$ and $i_0\ (1\eql i_0\eql n)$ be
an integer 
such that
\begin{eqnarray}
&&{\frac{{\rm Re}(a_{i_0})-1}{\alpha_{i_0}}}=\max_{1\eqls i
\eqls n}
\left[{\frac{{\rm Re}(a_i)-1}{\alpha_i}}\right].
\end{eqnarray}
Then there holds the asymptotic estimate$:$
\begin{eqnarray}
&&H^{m,n}_{p,q}(z)=h_{i_0}z^{(a_{i_0}-1)/\alpha_{i_0}}+o\left(z^{
(a_{i_0}-1)/
\alpha_{i_0}}\right)\quad(|z|\to\infty),
\end{eqnarray}
where $h_{i_0}$ is given by $(3.14)$ with $i=i_0$. In
particular},
\begin{eqnarray}
&&H^{m,n}_{p,q}(z)=O(z^{\varrho})\quad(|z|\to\infty)\quad\mbox{wi
th}\quad
\varrho=\max_{1\eqls i\eqls n}\left[{\frac{{\rm Re}(a_i)-1}
{\alpha_i}}\right].
\end{eqnarray}
\\\par
{\bf Remark 5.} \ The relations (3.4) and (3.12) are given in
[12, (2.2.4) 
and (2.2.7)] and [11, 8.3.2.3 and 8.3.2.4].
\\\\
\begin{flushleft}
{\bf 4. Explicit Power-Logarithmic Series Expansions}
\end{flushleft}
\setcounter{section}{4}
\setcounter{equation}{0}
\par
Now let us discuss the case when the condition (1.6) holds, but
(1.7) or 
(1.8) is violated:
\par
{\bf (e)} \ $\LL=\LL_{-\infty},\Delta\eqg0$ and some poles of the

Gamma-functions $\Gamma(b_j+\beta_js)\ (j=1,\cdots,m)$ coincide.
\par
{\bf (f)} \ $\LL=\LL_{+\infty},\Delta\eql0$  and some poles of
the 
Gamma-functions $\Gamma(1-a_i-\alpha_is)\ (i=1,\cdots,n)$
coincide. 
\par
First we consider the case (e). Let $b\equiv b_{jl}$ be one of
points (2.23) 
for which some poles of the Gamma-functions
$\Gamma(b_j+\beta_js)\ 
(j=1,\cdots,m)$ coincide and $N^*\equiv N^*_{jl}$ be order of
this pole. 
It means that there exist $j_1,\cdots,j_{N^*}\in \{1,\cdots,m\}$
and 
$l_{j_1},\cdots,l_{j_{N^*}}\in\N_0$ such that
\begin{eqnarray}
&&b=b_{jl}=-{\frac{b_{j_1}+l_{j_1}}{\beta_{j_1}}}=\cdots=
-{\frac{b_{j_{N^*}}+l_{j_{N^*}}}{\beta_{j_{N^*}}}}.
\end{eqnarray}
Then ${\HH}(s)z^{-s}$ has the pole of order $N^*$ at $b$ and
hence 
\begin{eqnarray}
&&\mathop{\rm
Res}_{s=b}[{\HH}(s)z^{-s}]={\frac1{(N^*-1)!}}\lim_{s\to b}
[(s-b)^{N^*}{\HH}(s)z^{-s}]^{(N^*-1)}.
\end{eqnarray}
\par
We denote
\begin{eqnarray}
&&{\HH}^*_1(s)=(s-b)^{N^*}\prod^{j_{N^*}}_{j=j_1}\Gamma(b_j+\beta
_js),
\qquad{\HH}^*_2(s)={\frac{{\HH}(s)}
{\displaystyle\prod^{j_{N^*}}_{j=j_1}\Gamma(b_j+\beta_js)}}.
\end{eqnarray}
Using the Leibniz rule, we have
\begin{eqnarray*}
&&[(s-b)^{N^*}{\HH}(s)z^{-s}]^{(N^*-1)}=\sum^{N^*-1}_{n=0}{N^*-1\
choose n}
[{\HH}^*_1(s)]^{(N^*-1-n)}[{\HH}^*_2(s)z^{-s}]^{(n)}\\[2mm]
&&\hspace{15mm}=\sum^{N^*-1}_{n=0}{N^*-1 \choose
n}[{\HH}^*_1(s)]^{(N^*-1-n)}
\sum^n_{i=0}{n \choose i}(-1)^i[{\HH}^*_2(s)]^{(n-i)}z^{-s}[\log
z]^i\\[2mm]
&&\hspace{15mm}=z^{-s}\sum^{N^*-1}_{i=0}\left\{\sum^{N^*-1}_{n=i}
(-1)^i
{N^*-1 \choose n}{n \choose
i}[{\HH}^*_1(s)]^{(N^*-1-n)}[{\HH}^*_2(s)]^{(n-i)}
\right\}[\log z]^i.
\end{eqnarray*}
Substituting this into (4.2), we obtain 
\begin{eqnarray}
&&\mathop{\rm Res}_{s=b_{jl}}[{\HH}(s)z^{-s}]=z^{(b_j+l)/\beta_j}
\sum^{N^*_{jl}-1}_{i=0}H^*_{jli}[\log z]^i,
\end{eqnarray}
where
\begin{eqnarray}
H^*_{jli}&\hspace{-2.5mm}\equiv&\hspace{-2.5mm}H^*_{jli}(N^*_{jl}
;b_{jl})
\\[2mm]
&\hspace{-2.5mm}=&\hspace{-2.5mm}{\frac1{(N^*_{jl}-1)!}}\sum^{N^*
_{jl}-1}_{n=i}
(-1)^i{N^*_{jl}-1\choose n}{n\choose
i}[{\HH}^*_1(b_{jl})]^{(N^*_{jl}-1-n)}
[{\HH}^*_2(b_{jl})]^{(n-i)}.\nonumber
\end{eqnarray}
\par
In particular, if $l=0$ and $i=N^*_{j0}-1$, then by setting 
$N^*_{j0}\equiv N^*_j$ and from (3.1), (4.1) and (4.3), we have
\begin{eqnarray}
H^*_j&\hspace{-2.5mm}\equiv &\hspace{-2.5mm}H^*_{j,0,N^*_j-1}(N^*_
j;b_{j0})=
{\frac{(-1)^{N^*_j-1}}{(N^*_j-1)!}}\
{\HH}^*_1(b_{j0}){\HH}^*_2(b_{j0})\\[2mm]
&\hspace{-2.5mm}=&\hspace{-2.5mm}{\frac{(-1)^{N^*_j-1}}{(N^*_j-1)
!}}
\left\{\prod^{N^*_j}_{k=1}{\frac{(-1)^{j_k}}{{j_k}!\beta_{j_k}}}
\right\}
{\frac{\displaystyle{\prod^m_{i=1,i\neq j_1,\cdots,j_{N^*_j}}
\Gamma\left(b_i-\frac{b_j\beta_i}{\beta_j}\right)\prod^n_{i=1}
\Gamma\left(1-a_i+\frac{b_j\alpha_i}{\beta_j}\right)}}
{\displaystyle{\prod^p_{i=n+1}\Gamma\left(a_i-\frac{b_j\alpha_i}
{\beta_j}
\right)\prod^q_{i=m+1}\Gamma\left(1-b_i+\frac{b_j\beta_i}{\beta_j}
\right)}}}.
\nonumber
\end{eqnarray}
\par
Thus, in view of Theorem 2(A), we have
\\\par
{\bf Theorem 5.} \ {\sl Let the conditions in $(1.6)$ be
satisfied and let 
either $\Delta>0,z\neq0$ or $\Delta=0,0<|z|<\delta$. Then the
$H$-function 
$(1.1)$ has the power-logarithmic series expansion   
\begin{eqnarray}
&&H^{m,n}_{p,q}(z)={\sum_{j,l}}'\ h^*_{jl}z^{(b_j+l)/\beta_j}+
{\sum_{j,l}}''\sum^{N^*_{jl}-1}_{i=0}H^*_{jli}z^{(b_j+l)/\beta_j}
[\log z]^i.
\end{eqnarray}
Here $\sum'$ and $\sum''$ are summations taken over 
$j,l\ (j=1,\cdots,m; l=0,1,\cdots )$ such that the
Gamma-functions 
$\Gamma(b_j+\beta_js)$ have simple poles and poles of order
$N^*_{jl}$ at the 
points $b_{jl},$ respectively$,$ and the constants $h^*_{jl}$ are
given by 
$(3.3)$ while the constants $H^*_{jli}$ are given by} (4.5).
\\\par
{\bf Corollary 6.} \ {\sl If the conditions in $(1.6)$ are
satisfied and 
$\Delta\eqg0,$ then $(4.7)$ gives the asymptotic expansion of 
$H^{m,n}_{p,q}(z)$ near zero and the main terms of this
asymptotic formula 
have the form$:$
\begin{eqnarray}
H^{m,n}_{p,q}(z)&\hspace{-2.5mm}=&\hspace{-2.5mm}\mathop{\mbox{$
{\displaystyle\sum}'$}}^m_{j=1}\left[h^*_jz^{b_j/\beta_j}+O(z^{(b
_j+1)/
\beta_j})\right]\\[2mm]
&&+\mathop{\mbox{${\displaystyle\sum}''$}}^{m}_{j=1}\left(H^*_jz^
{b_j/
\beta_j}[\log z]^{N^*_j-1}+O(z^{(b_j+1)/\beta_j}[\log
z]^{N^*_j-1})\right)
\quad(z\to0).\nonumber
\end{eqnarray}
Here $\sum'$ and $\sum''$ are summations taken over $j\
(j=1,\cdots,m)$ such 
that the Gamma-functions $\Gamma(b_j+\beta_js)$ have simple poles
and poles 
of order $N^*_j\equiv N^*_{j0}$ at the points $b_{j0},$
respectively$,$ and 
$h^*_{j}$ are given by $(3.6)$ while $H^*_{j}$ are given by}
(4.6).
\\\par
{\bf Corollary 7.} \ {\sl Let the conditions in $(1.6)$ be
satisfied$,$ and 
let $\Delta\eqg0$ and $b_{jl}$ be poles of the Gamma-functions 
$\Gamma(b_j+\beta_js)\ (j=1,\cdots,m).$ Let $j_{01}$ and 
$j_{02}\ (1\eql j_{01},j_{02}\eql m)$ be integers such that
\begin{eqnarray}
&&{\frac{{\rm Re}(b_{j_{01}})}{\beta_{j_{01}}}}=\min_{1\eqls
j\eqls m}
\left[{\frac{{\rm Re}(b_{j})}{\beta_{j}}}\right],
\end{eqnarray}
when the poles $b_{jl}\ (j=1,\cdots ,m;l\in\N_0)$ are simple, and

\begin{eqnarray}
&&{\frac{{\rm Re}(b_{j_{02}})}{\beta_{j_{02}}}}=\min_{1\eqls
j\eqls m}
\left[{\frac{{\rm Re}(b_{j})}{\beta_{j}}}\right],
\end{eqnarray}
when the poles $b_{jl}\ (j=1,\cdots ,m;l\in\N_0)$ coincide. 
\par
{\bf a)} \ If $j_{01}<j_{02},$ then the asymptotic expansion of
the 
$H$-function has the form$:$
\begin{eqnarray}
&&H^{m,n}_{p,q}(z)=h^*_{j_{01}}z^{b_{j_{01}}/\beta_{j_{01}}}
+o\left(z^{b_{j_{01}}/\beta_{j_{01}}}\right)\quad(z\to 0),
\end{eqnarray}
where $h^*_{j_{0}}$ is given by $(3.6)$ with $j=j_{01}$. In
particular$,$ 
the relation $(3.9)$ holds.
\par
{\bf b)} \ If $j_{01}\eqg j_{02}$ and $b_{j_{02},0}$ has the pole
of order 
$N^*_{j_{02}},$ then the first term in asymptotic expansion of
the 
$H$-function has the form$:$
\begin{eqnarray}
&&H^{m,n}_{p,q}(z)=H^*_{j_{02}}z^{b_{j_{02}}/\beta_{j_{02}}}
[\log (z)]^{N^*_{j_{02}}-1}+o\left(z^{b_{j_{02}}/\beta_{j_{02}}}
[\log (z)]^{N^*_{j_{02}}-1}\right)\quad(z\to 0),
\end{eqnarray}
where $H^*_{j_{02}}$ is given by $(4.6)$ with $j=j_{02}$. In
particular$,$ if 
$N^*$ is the largest order of general poles of Gamma-functions 
$\Gamma(b_j+\beta_js)\ (j=1,\cdots,m),$ then}
\begin{eqnarray}
&&H^{m,n}_{p,q}(z)=O\left(z^{\rho}[\log
(z)]^{N^*}\right)\quad(z\to 0)
\quad\mbox{with}\quad\rho=\min_{1\eqls j\eqls m}
\left[{\frac{{\rm Re}(b_{j})}{\beta_{j}}}\right].
\end{eqnarray}
\\\par
Now we consider the case (f). Let $a=a_{ik}$ be one of points
(2.24) for 
which some poles of $\Gamma(1-a_i-\alpha_is)$ coincide and 
$N=N_{ik}$ be order of this pole. It means that there exist 
$i_{1},\cdots,i_{N}\in \{1,\cdots ,n\}$ and
$k_{i_{1}},\cdots,k_{i_{N}}
\in\N_0$ such that
\begin{eqnarray}
&&a=a_{ik}\equiv
{\frac{1-a_{i_{1}}+k_{i_{1}}}{\alpha_{i_{1}}}}=\cdots
={\frac{1-a_{i_{N}}+k_{i_{N}}}{\alpha_{i_{N}}}}.
\end{eqnarray}
Then the integrand ${\HH}(s)z^{-s}$ of the integral (1.1) has the
pole of 
order $N$ at $a$. Similarly to (4.3), we denote
\begin{eqnarray}
&&{\HH}_{1}(s)=(s-a)^{N}\prod^{i_{N}}_{i=i_{1}}\Gamma(1-a_i-
\alpha_i s),
\ {\HH}_{2}(s)={\frac{{\HH}(s)}
{\displaystyle{\prod^{i_{N}}_{i=i_{1}}\Gamma(1-a_i-\alpha_is)}}}
\end{eqnarray}
and, then, find similarly to (4.4) and (4.5) that 
\begin{eqnarray}
&&\mathop{\rm
Res}_{s=a_{ik}}[{\HH}(s)z^{-s}]=z^{(a_{i}-1-k)/\alpha_{i}}
\sum^{N_{ik}-1}_{j=0}H_{ikj}[\log (z)]^{j},
\end{eqnarray}
where
\setcounter{equation}{16}
\begin{eqnarray}
H_{ikj}&\hspace{-2.5mm}\equiv&\hspace{-2.5mm}H_{ikj}(N_{ik};a_{ik
})\\[2mm]
&\hspace{-2.5mm}=&\hspace{-2.5mm}{\frac{1}{(N_{ik}-1)!}}
\left\{\sum^{N_{ik}-1}_{n=j}(-1)^{j}{N_{ik}-1 \choose n}{n
\choose j}
[{\HH}_1(a_{ik})]^{(N_{ik}-1-n)}[{\HH}_2(a_{ik})]^{(n-j)}\right\}
.\nonumber
\end{eqnarray}
\par
In particular, if we set $k=0$, $j=N_{i0}-1$ and $N_{i0}\equiv
N_i$, then, 
using (4.15) and (3.1), we have
\begin{eqnarray}
&&\hspace{-3mm}H_i\equiv
H_{i,0,N_i-1}(N_i;a_{i0})={\frac{(-1)^{N_i-1}}
{(N_i-1)!}}{\HH}_1(a_{i0}){\HH}_2(a_{i0})={\frac{(-1)^{N_i-1}}
{(N_i-1)!}}\\[2mm]
&&\times
\left(\prod^{N_i}_{k=1}{\frac{(-1)^{{i_k}-1}}{{i_k}!\alpha_
{i_k}}}
\right)
{\frac{\displaystyle{\prod^m_{j=1}\Gamma\left(b_j+[1-a_i]
\frac{\beta_j}
{\alpha_i}\right)
\prod^n_{j=1,j\neq
i_1,\cdots,i_{N_i}}\Gamma\left(1-a_j-[1-a_i]\frac{\alpha_j}
{\alpha_i}\right)}}
{\displaystyle{\prod^p_{j=n+1}\Gamma\left(a_j+[1-a_i]
\frac{\alpha_j}{\alpha_i}
\right)
\prod^q_{j=m+1}\Gamma\left(1-b_j-[1-a_i]\frac{\beta_j}{\alpha_i}
\right)}}}.
\nonumber
\end{eqnarray}
\par
Therefore, Theorem 2(B) implies the similar result to Theorem 5:
\\\par
{\bf Theorem 6.} \ {\sl Let the conditions in $(1.6)$ be
satisfied and let 
either $\Delta<0,z\neq0$ or $\Delta=0,|z|>\delta$. Then the
$H$-function 
$(1.1)$ has the power-logarithmic series expansion   
\begin{eqnarray}
&&H^{m,n}_{p,q}(z)={\sum_{i,k}}'\
h_{ik}z^{(a_{i}-1-k)/\alpha_{i}}
+{\sum_{i,k}}''\ \sum^{N_{ik}-1}_{j=0}H_{ikj}z^{(a_{i}-1-k)/
\alpha_{i}}[\log (z)]^{j}.
\end{eqnarray}
Here $\sum'$ and $\sum''$ are summations taken over 
$i,k\ (i=1,\cdots,n;\ k=0,1,\cdots )$ such that Gamma-functions 
$\Gamma(1-a_i-\alpha_is)$ have simple poles and poles of order
$N_{ik}$ at 
the points $a_{ik},$ respectively$,$ and the constants $h_{ik}$
are given by 
$(3.11)$ while the constants $H_{ikj}$ are given by} (4.17).
\\\par
{\bf Corollary 8.} \ {\sl If the conditions in $(1.6)$ are
satisfied and 
$\Delta\eql0,$ then $(4.19)$ gives the asymptotic expansion of 
$H^{m,n}_{p,q}(z)$ near infinity and the main terms of this
asymptotic formula 
have the form$:$
\setcounter{equation}{19}
\begin{eqnarray}
&&\hspace{-10mm}H^{m,n}_{p,q}(z)=\mathop{\mbox{${\displaystyle
\sum}'$}}^n_{i=1}
\left[h_iz^{(a_i-1)/\alpha_i}+O(z^{(a_i-2)/\alpha_i})
\right]\\[2mm]
&&+\mathop{\mbox{${\displaystyle\sum}''$}}^n_{i=1}
\left(H_iz^{(a_i-1)/\alpha_i}
[\log(z)]^{N_i-1}+O(z^{(a_i-2)/\alpha_i}
[\log z]^{N_i-1})\right)\quad(|z|\to\infty).\nonumber
\end{eqnarray}
Here $\sum'$ and $\sum''$ are summations taken over $i\
(i=1,\cdots,n)$ such 
that the Gamma-functions $\Gamma(1-a_i-\alpha_is)$ have simple
poles and poles 
of order $N_i\equiv N_{i0}$ in the points $a_{i0}$ in $(2.24),$ 
respectively$,$ and $h_{i}$ are given by $(3.14)$ while $H_{i}$
are given by} 
(4.18).
\\\par
{\bf Corollary 9.} \ {\sl Let the conditions in $(1.6)$ be
satisfied$,$ and 
let $\Delta\eql0$ and $a_{ik}$ be poles of the Gamma-function 
$\Gamma(1-a_i-\alpha_is)\ (i=1,\cdots,n).$ Let $i_{01}$ and 
$i_{02}\ (1\eql i_{01},i_{02}\eql n)$ be integers such that
\begin{eqnarray}
&&{\frac{{\rm Re}(a_{i_{01}})-1}{\alpha_{i_{01}}}}=\max_{1\eqls
i\eqls n}
\left[{\frac{{\rm Re}(a_i)-1}{\alpha_i}}\right],
\end{eqnarray}
when the poles $a_{ik}\ (i=1,\cdots,n;k=0,1,\cdots)$ are
simple$,$ and 
\begin{eqnarray}
&&{\frac{{\rm Re}(a_{i_{02}})-1}{\alpha_{i_{02}}}}=\max_{1\eqls
i\eqls n}
\left[{\frac{{\rm Re}(a_i)-1}{\alpha_i}}\right],
\end{eqnarray}
when the poles $a_{ik}\ (i=1,\cdots,n;k\in\N_0)$ coincide. 
\par
{\bf (a)} \ If $i_{01}<i_{02}$, then the asymptotic expansion of
the 
$H$-function has the form$:$
\begin{eqnarray}
&&H^{m,n}_{p,q}(z)=h_{i_{01}}z^{(a_{i_{01}}-1)/\alpha_{i_{01}}}
+o\left(z^{(a_{i_{01}}-1)/\alpha_{i_{01}}}\right)\quad(|z|\to\inf
ty),
\end{eqnarray}
where $h_{i_{01}}$ is given by $(3.14)$ with $i=i_{01}$. In
particular$,$ 
the relation $(3.17)$ holds.
\par
{\bf (b)} \ If $i_{01}\eqg i_{02}$ and $a_{i_{02},0}$ has the
pole of order 
$N_{i_{02}},$ then the asymptotic expansion of the $H$-function
has the 
form$:$
\begin{eqnarray}
&&H^{m,n}_{p,q}(z)=H_{i_{02}}z^{(a_{i_{02}}-1)/\alpha_{i_{02}}}
[\log z]^{N_{i_{02}}-1}+o\left(z^{(a_{i_{02}}-1)/\alpha_{i_{02}}}
[\log (z)]^{N_{i_{02}}-1}\right)\quad(|z|\to\infty),
\end{eqnarray}
where $H_{i_{02}}$ is given by $(4.18)$ with $i=i_{02}$. In
particular$,$ if 
$N$ is the smallest order of general poles of Gamma-functions 
$\Gamma(1-a_i-\alpha_is)\ (i=1,\cdots,n),$ then}
\begin{eqnarray}
&&H^{m,n}_{p,q}(z)=O\left(z^{\varrho}[\log
(z)]^{N}\right)\quad(|z|\to\infty)
\quad\mbox{with}\quad\varrho=\min_{1\eqls i\eqls
n}\left[{\frac{{\rm Re}
(a_{i})-1}{\alpha_{i}}}\right].
\end{eqnarray}
\\\par
In conclusion, we give the following consequence of Corollaries
3, 5, 7 and 
9, which unifies the power and power-logarithmic asymptotic
behavior of 
$H^{m,n}_{p,q}(z)$ near zero and infinity.
\\\par
{\bf Theorem 7.} \ {\sl Let the conditions in $(1.6)$ be
satisfied.
\par
{\bf (a)} \ If $\Delta\eqg0$ and the poles of the Gamma-function 
$\Gamma(b_j+\beta_js)\ (j=1,\cdots,m)$ are simple$,$ then the
$H$-function 
$(1.3)$ has the asymptotic estimate $(3.9)$ at zero. If some of
poles of 
$\Gamma(b_j+\beta_js)\ (j=1,\cdots,m)$ coincide$,$ then 
$H^{m,n}_{\thinspace p,q}(z)$ has the asymptotic estimate either
$(3.9)$ or 
$(4.13)$ at zero.
\par
{\bf (b)} \ If $\Delta\eql0$ and the poles of the Gamma-function 
$\Gamma(1-a_i-\alpha_is)\ (i=1,\cdots,n)$ are simple$,$ then the 
$H$-function $(1.3)$ has the asymptotic estimate $(3.17)$ at
infinity. If 
some of poles of $\Gamma(1-a_i-\alpha_is)\ (i=1,\cdots,n)$
coincide$,$ then 
$H^{m,n}_{\thinspace p,q}(z)$ has the asymptotic estimate
$(3.17)$ or 
$(4.25)$ at infinity.}
\\\par
{\bf Remark 6.} \ The power-logarithmic expansions and more
complicated 
results than in (4.7) were indicated in [9, Section\thinspace3.7]
(see also 
[8, Section\thinspace5.8]) and the particular cases
$H^{p,0}_{0,p}(z)$ and 
$H^{p,0}_{p,p}(z)$ in [7].
\\\\
\begin{flushleft}
{\bf Acknowledgement}
\end{flushleft}
\par
The authors would like to express their gratitude to the referee
for valuable 
comments.
\vspace{1cm}
\begin{center}{\bf References}\end{center}
\baselineskip=4.47mm
\begin{enumerate}
\item[{[1]\ }] Braaksma, B.L.G., Asymptotic expansions and
analytic 
continuation for a class of Barnes-integrals, {\it Compositio
Math.} 
{\bf 15}(1964), 239-341.
\item[{[2]\ }] Dixon, A.L. and Ferrar, W.L., A class of
discontinuous 
integrals, {\it Quart. J. Math.$,$ Oxford Ser.} {\bf 7}(1936),
81-96.
\item[{[3]\ }] Erd\'elyi, A., Magnus, W., Oberhettinger, F. and
Tricomi, 
F.G., {\it Higher Transcendental Functions,\ Vol.\thinspace I}, 
McGraw-Hill, New York, Toronto and London, 1953.    
\item[{[4]\ }] Fox, C., The $G$ and $H$ functions as symmetrical
Fourier 
kernels., {\it Trans. Amer. Math. Soc.} {\bf 98}(1961), 395-429.
\item[{[5]\ }] Kilbas, A.A., Saigo, M. and Shlapakov, S.A.,
Integral 
transforms with Fox's $H$-function in spaces of summable
functions, 
{\it Integral Transform. Spec. Funct.} {\bf 1}(1993), 87-103.
\item[{[6]\ }] Marichev, O.I., {\it Handbook of Integral
Transforms of 
Higher Transcendental Function. Theory and Algorithmic Tables}, 
Wiley (Ellis Horwood), New York, Brisbane, Chichester and
Toronto, 1982.
\item[{[7]\ }] Mathai, A.M., An expansion of Meijer's
$G$-function in the 
logarithmic case with applications, {\it Math. Nachr.} {\bf
48}(1971), 129-139.
\item[{[8]\ }] Mathai, A.M., A few results on the exact
distributions of 
certain multivariete statistics. II, {\it Multivariate
Statistical Inference} 
(Proc. Res. Sem. Dalhouse Univ., Halifax, N.S., 1972), 169-181,
North-Holland, 
Amsterdam and New York, 1973.
\item[{[9]\ }] Mathai, A.M. and Saxena, R.K., {\it The
$H$-Function with 
Applications in Statistics and other Disciplines}, Wiley (Halsted
P.), New 
York, London, Sydney and Toronto, 1978.
\item[{[10]\ }] Mellin, Hj., Abri\ss\ einer einheitlichen Theorie
der Gamma- 
und der hypergeometrischen Funktionen, {\it Math. Ann.} {\bf
68}(1910), 305-337.
\item[{[11]\ }] Prudnikov, A.P., Brychkov, Yu.A. and Marichev,
O.I., 
{\it Integrals and Series, Vol.3, More Special Functions}, Gordon
and Breach, 
New York, Philadelphia, London, Paris, Montreux, Tokyo and
Melbourne, 1990.  
\item[{[12]\ }] Srivastava, H.M., Gupta, K.C. and Goyal, S.P.,
{\it The 
$H$-Functions of One and Two Variables with Applications}, South
Asian 
Publishers, New Delhi and Madras, 1982.
\end{enumerate}

\end{document}